\documentclass[12pt]{book}
\usepackage{makeidx}
\usepackage{multirow}
\usepackage{multicol}
\usepackage[dvipsnames,svgnames,table]{xcolor}
\usepackage{graphicx}
\usepackage{epstopdf}
\usepackage{ulem}
\usepackage{tikz}
\usepackage{fontenc}
\usepackage{mathrsfs}
\usepackage{float}
\usepackage{curves}
\usepackage{caption}
\usepackage{hyperref}
\usepackage{amsmath,amsfonts,amssymb,amscd,amsmath,enumerate,verbatim}
\usepackage{amssymb}

\usepackage[paperwidth=560pt,paperheight=700pt,top=72pt,right=67pt,bottom=72pt,left=72pt]{geometry}
\begin{document}
\thispagestyle{empty}
\setcounter{page}{01}
\def\ord{\rm ord}
\begin{large}
\begin{center}
\textbf{A note on Apery's constant is transcendental}
\end{center}

\begin{center}
{\bf Shekhar Suman}\vspace{.2cm}\\
M.Sc., I.I.T. Delhi, India\vspace{.2cm}\\
Email: shekharsuman068@gmail.com\\

\end{center}\vspace{.1cm}
{\bf Abstract:} Beuker's [2] considers the following integral $$
\int_{0}^{1}\int_{0}^{1} \frac{-\log xy}{1-xy} P_n(x)P_n(y)\ dx dy$$If $d_n=\text{LCM}(1,2,...,n)$, then $$
    0<\frac{|A_n+B_n\zeta(3)|}{d_n^3}<2(\sqrt{2}-1)^{4n} \zeta(3)
$$ for some $A_n,B_n\in\mathbb{Z}$. We establish that if Apery's constant is algebraic then the above inequality fails to be true. This proves that $\zeta(3)$ is transcendental.  
\vspace{0.2cm}\vspace{.2cm} \vspace{0.2cm}\\
\noindent{\bf 1. Introduction and Definitions} The Riemann zeta function, $\zeta(s)$ is defined as the analytic continuation  of the Dirichlet series
$$\zeta(s) = \sum_{n=1}^{{\infty}} \frac{1}{n^s}$$ which converges in the half plane $\Re(s)>1$.\ The \ Riemann \ zeta \ function \ is \ a \ meromorphic \ function \ on \ the \ whole \ complex \ plane, \ which \ is \ holomorphic \ everywhere \ except \ for \ a \ simple \ pole \ at \ s = 1 \ with \ residue \ 1.\ The Riemann Hypothesis \ states \ that \ all \ the \ non \ trivial \ zeros \ of \ the \ Riemann \ zeta \ function \ lie \ on \ the \ critical \ line \ $\Re(s)=\frac{1}{2}$. The Riemann xi function is defined as $$ \xi(s)=\frac{1}{2}s(s-1)\pi^{-s/2}\Gamma\left(\frac{s}{2}\right)\zeta(s)$$ $\xi(s)$ is an entire function whose zeros are the non trivial zeros of $\zeta(s)$ [1]. Further $\xi(s)$ satisfies the functional equation [1] $$\xi(s)=\xi(1-s)$$ It is known that the Apery's constant, $\zeta(3)$ is irrational [2].\\\\
{\bf 2. Main Theorems} 
The goal of this article is to prove the following result.\\\\
\textbf{\underline{Theorem :}} Apery's constant, $\zeta(3)$ is transcendental.\\\\
\textbf{Proof}: We prove the claim by the method of contradiction. Let us assume that $\zeta(3)$ is algebraic. Then it satisfies \begin{equation}
    c_0+c_1\zeta(3)+c_2\zeta^2(3)+...+c_m\zeta^m(3) =0
\end{equation} where $c_i\in\mathbb{Z}$, $c_0>0$ and $m\in\mathbb{Z}^+$.\\ It is well known that [2] \begin{equation}
    2\zeta(3)=\int_{0}^{1}\int_{0}^{1} \frac{-\log xy}{1-xy} \ dx dy
\end{equation}
Beuker's [2] considers the following integral \begin{equation}
    I_n:= \int_{0}^{1}\int_{0}^{1} \frac{-\log xy}{1-xy} P_n(x)P_n(y)\ dx dy
\end{equation}
By Beuker's [2] proof, we have if $d_n=\text{LCM}(1,2,...,n)$, then \begin{equation}
    0<\frac{|A_n+B_n\zeta(3)|}{d_n^3}<2(\sqrt{2}-1)^{4n} \zeta(3)
\end{equation} for some $A_n,B_n\in\mathbb{Z}$. Equation (2) can be rewritten as \begin{equation}
    0<\frac{|A_n+B_n\zeta(3)|}{2(\sqrt{2}-1)^{4n} d_n^3}<\zeta(3)
\end{equation} 
So we get on squaring equation (3) \begin{equation}
    0<\frac{|A_n+B_n\zeta(3)|^2}{(2(\sqrt{2}-1)^{4n} d_n^3)^2}<\zeta^2(3)
\end{equation} Similarly we get \begin{equation}
    0<\frac{|A_n+B_n\zeta(3)|^3}{(2(\sqrt{2}-1)^{4n} d_n^3)^3}<\zeta^3(3)
\end{equation}and so on \begin{equation*}
    0<\frac{|A_n+B_n\zeta(3)|^{m+1}}{(2(\sqrt{2}-1)^{4n} d_n^3)^{m+1}}<\zeta^{m+1}(3)\tag{$m$+5}
\end{equation*} Now multiplying (3) by $c_0$, (4) by $c_1$, (5) by $c_2$,..., ($m$+3) by $c_m$, we get \begin{align*}
    0&<c_0\frac{|A_n+B_n\zeta(3)|}{2(\sqrt{2}-1)^{4n} d_n^3}+c_1\frac{|A_n+B_n\zeta(3)|^2}{(2(\sqrt{2}-1)^{4n} d_n^3)^2}+c_2 \frac{|A_n+B_n\zeta(3)|^3}{(2(\sqrt{2}-1)^{4n} d_n^3)^3}\notag\\&
    +...+c_m\frac{|A_n+B_n\zeta(3)|^{m+1}}{(2(\sqrt{2}-1)^{4n} d_n^3)^{m+1}}<c_0\zeta(3)+c_1\zeta^2(3)+c_2\zeta^3(3)+...+c_m\zeta^{m+1}(3)
\tag{$m$+6}\end{align*} 
\begin{align*}
    0&<c_0\frac{|A_n+B_n\zeta(3)|}{2(\sqrt{2}-1)^{4n} d_n^3}+c_1\frac{|A_n+B_n\zeta(3)|^2}{(2(\sqrt{2}-1)^{4n} d_n^3)^2}+c_2 \frac{|A_n+B_n\zeta(3)|^3}{(2(\sqrt{2}-1)^{4n} d_n^3)^3}\notag\\&
    +...+c_m\frac{|A_n+B_n\zeta(3)|^{m+1}}{(2(\sqrt{2}-1)^{4n} d_n^3)^{m+1}}<\zeta(3)(c_0+c_1\zeta(3)+c_2\zeta^2(3)+...+c_m\zeta^{m}(3))
\tag{$m$+7}\end{align*} 
Now using equation (1) we get \begin{align*}
    0&<c_0\frac{|A_n+B_n\zeta(3)|}{2(\sqrt{2}-1)^{4n} d_n^3}+c_1\frac{|A_n+B_n\zeta(3)|^2}{(2(\sqrt{2}-1)^{4n} d_n^3)^2}+c_2 \frac{|A_n+B_n\zeta(3)|^3}{(2(\sqrt{2}-1)^{4n} d_n^3)^3}\notag\\&
    +...+c_m\frac{|A_n+B_n\zeta(3)|^{m+1}}{(2(\sqrt{2}-1)^{4n} d_n^3)^{m+1}}<0
\tag{$m$+8}\end{align*} So we get \begin{equation*}
    0<\sum_{k=1}^{m+1} c_{k-1}\frac{|A_n+B_n\zeta(3)|^k}{(2(\sqrt{2}-1)^{4n} d_n^3)^{k}}<0
\tag{$m$+9}\end{equation*}
which is a contradiction. So our assumption that $\zeta(3)$ is algebraic is not possible. Hence $\zeta(3)$ is transcendental. 
This settles the proof of the Theorem.\\\\
\\\\
\\\\
\\\\
\\\\
\\\\
\\\\
\\\\
\\\\
\begin{center}
{\bf References}
\end{center} 
\big[1\big] Edwards, H.M., {\it Riemann's\  Zeta \ Function} , Dover Publications (2001).\\\\
\big[2\big] Beukers, F., {\it A note on the irrationality of $\zeta(2)$ and $\zeta(3)$}, Bulletin of the London Mathematical Society, Vol. 11, issue 3, 268-272 (1979).\\\\
\big[3\big] Apéry, R. (1979). {\it Irrationalité \ de \ $\zeta(2)$ \ et \ $\zeta(3)$}, Astérisque. 61: 11–13.
\\\\ \vspace{0.2cm}\\\\
 \noindent{\bf Data Availability Statement (DAS)}\\\\
\big[1\big] Books available in the library or online library - for the reference [1](\href{http://libgen.rs/book/index.php?md5=2EBBD146008A02EB5F22CD97F3E4DC37}{Click here for [1]})\\\\
\big[2\big] Data openly available in a public repository that does not issue DOIs - for the reference, [2] (\href{https://londmathsoc.onlinelibrary.wiley.com/doi/abs/10.1112/blms/11.3.268}{Click here for [2]}) and [3] (\href{https://www.scirp.org/(S(czeh2tfqw2orz553k1w0r45))/reference/referencespapers.aspx?referenceid=2972936}{Click here for [3]}) \end{large}
\end{document}